\documentclass[10pt,onecolumn,draftclsnofoot]{IEEEtran}
\usepackage[utf8]{inputenc}

\usepackage{lipsum}
\usepackage{amsfonts, amssymb, amsmath, mathtools}
\usepackage{graphicx}
\usepackage{algorithmic}
\usepackage{xcolor}

\usepackage{enumitem}
\usepackage{subcaption}
\usepackage[linesnumbered,algoruled,boxed,lined]{algorithm2e}

\usepackage{authblk}

\DeclareMathOperator*{\argmin}{arg\,min}

\newcommand{\mb}[1]{\mathbf{#1}}

\usepackage{enumitem}
\setlist[enumerate]{leftmargin=.5in}
\setlist[itemize]{leftmargin=.5in}

\newtheorem{thm}{Theorem}

\newtheorem{definition}[thm]{Definition}
\newtheorem{lemma}[thm]{Lemma}

\title{Efficient Projection onto the $\ell_{\infty,1}$ Mixed-Norm Ball using a Newton Root Search Method
}

\author[1]{Gustavo Chau}
\author[2]{Brendt Wohlberg}
\author[1]{Paul Rodriguez}
\affil[1]{Electrical Engineering Department, Pontificia Universidad Cat\'olica del Per\'u, Lima, Per\'u
\authorcr Email: {\tt \{gustavo.chau, prodrig\}@pucp.edu.pe}\vspace{1.5ex}}
\affil[2]{Theoretical Division, Los Alamos National Laboratory, Los Alamos, NM 87545 USA 
\authorcr Email: {\tt brendt@lanl.gov} \vspace{-2ex}}

\usepackage{amsopn}

\usepackage{amssymb}
\usepackage{tabularx}

 \newcolumntype{Y}{>{\centering\arraybackslash}X}
 \newcolumntype{R}{>{\raggedleft\arraybackslash}X}
 \newcolumntype{M}[1]{>{\centering\arraybackslash}m{#1}}
 \usepackage{color}

\usepackage[font=small,skip=4pt]{caption}

\DeclareMathOperator{\shrink}{shrink}

\begin{document}


%
\maketitle

\begin{abstract}

Mixed norms that promote structured sparsity have numerous applications in signal processing and machine learning problems. In this work, we present a new algorithm,  based on a Newton root search technique, for computing the projection onto the $\ell_{\infty,1}$ ball, which has found application in cognitive neuroscience and classification tasks. Numerical simulations show that our proposed method is between 8 and 10 times faster on average, and up to 20 times faster for very sparse solutions, than the previous state of the art.
Tests on real functional magnetic resonance image data show that, for some data distributions, our algorithm can obtain speed improvements by a factor {of between 10 and 100, depending on the implementation}. 
\end{abstract}


\section{Introduction}

Mixed norms are important in modeling group correlations in applications such as genetics~\cite{yuan2011efficient}, 
electroencephalography~\cite{gramfort2012mixed} and signal processing~\cite{kowalski2009sparse}. In this work, 
we consider mixed norms with non-overlapping groups applied to matrix-form data $A \in \mathbb{R}^{M \times N}$, 
where the rows $\mathbf{a}_m \in \mathbb{R}^{N}$ represent the different groups. 
Following the notation of \cite{kowalski2009sparse}, we define the $\ell_{p,q}$-norm of $A$ as 
\begin{equation}
	\| A \|_{p,q} = \left(\sum_{m=1}^{M}{||\mathbf{a}_m||_p^q}\right)^{1/q}. 
\end{equation}
We will focus on a special case, the $\ell_{\infty,1}$-norm:
  \begin{equation}
      \| A \|_{\infty,1} = \sum_{m=1}^{M}{\| \mathbf{a}_m \|_{\infty}} \;,
      \label{eqn:infty-1}
  \end{equation}
where $\| \mathbf{u} \|_{\infty} = \max_n \{ |u_n| \}$ for $\mathbf{u} \in \mathbb{R}^{N}$.

The main contribution of this work is a new, computationally efficient algorithm for computing the projection 
onto the $\ell_{\infty,1}$-ball:
\begin{align}
\mathrm{proj}_{\|\cdot\|_{\infty,1}}(B,\tau) & := 
\argmin_X \; \frac{1}{2} \| X - B  \|_F^2 \label{eq:proj_1inf_ball1}\\ 
\text{ s.t. } &  \| X \|_{\infty,1} \leq \tau \;. \notag
\end{align}

This $\ell_{\infty,1}$ constraint problem has been applied to image annotation~\cite{quattoni2009efficient}, cognitive neuroscience~\cite{liu2009blockwise} and least absolute shrinkage and selection operator operator (LASSO) regression~\cite{sra2011fast}.
We propose a novel approach for solving~\eqref{eq:proj_1inf_ball1} that utilizes a root search based on a Newton method, in which the total number of major iterations for the root search is reduced by applying a simple scheme for choosing a feasible initial solution. This approach provides significantly improved performance compared with our previous approach based on a Steffensen root finding method~\cite{chau2018steffensen}. 
{The contributions of the present manuscripts can be summarized as:}
\begin{enumerate}[label=(\roman*)]
\item We present a new method for solving \eqref{eq:proj_1inf_ball1}. Instead of a Steffensen root-finding procedure, we formally develop an approximated Newton method for the root search function. 
\item We significantly expand the theoretical analysis of the initial point estimation and pruning.
\item We consider additional experiments conducted on real functional magnetic resonance imaging (fMRI) data in order to validate the usefulness of our proposed method.
\end{enumerate}
The manuscripts is organized as follows: Section \ref{sec:existing} summarizes existing approaches for solving~\eqref{eq:proj_1inf_ball1}, Section \ref{sec:preliminaries} presents some mathematical preliminaries needed for our derivations, Section \ref{sec:proposed} describes our proposed method, Sections \ref{sec:mtl} and \ref{sec:rslts} present our results on simulated and real data, and in Sections \ref{sec:discussion} and \ref{sec:conclusions}, we discuss these results and present the conclusions of our work.

\section{Existing approaches}
\label{sec:existing}
In this section, we review some previous approaches for solving~\eqref{eq:proj_1inf_ball1}. These methods tend to focus on
reinterpreting the problem via some form {of} quadratic or linear programming (Sections~\ref{sec:turlach} and 
\ref{sec:quattoni}) or via root-finding (Sections~\ref{sec:sra} and \ref{sec:steff}).

\subsection{Solution via interior point methods}
\label{sec:turlach}
{An approach for achieving simultaneous variable selection by considering the problem}
 \begin{align}
    \min_X \; &
    \frac{1}{2} \| QX - B  \|_F^2 \label{eq:proj_1inf_ball1_turlach}\\ 
\text{ s.t. } &  \| X \|_{\infty,1} \leq \tau \;, \notag
 \end{align}
 where $Q$ is a fixed matrix and {$\tau$ is a problem parameter, was proposed in \cite{turlach2005simultaneous} }. For the particular case where $Q=I$, this approach consisted of introducing the variables $\rho_m, m \in \{1,\ldots,M \}$ and recasting the problem as a convex quadratic optimization:
 \begin{align}
    \min_{\{ \rho_m \}} \; &
    \frac{1}{2} \sum_{n=1}^N \sum_{m=1}^M (|b_{nm}|-\rho_m)^2_+ \label{eq:proj_1inf_ball1_turlach2}\\ 
\text{ s.t. } & \sum_{m=1}^M \rho_m = \tau \;,\notag \\
	& \rho_m \geq 0 \;\;, \quad m \in \{1,\ldots,M \} \;, \notag
 \end{align}
where $(x)_+ = \max(x,0)$. It was shown that the $\rho_l$ are piecewise linear functions of $\tau$ and that they fulfill the Karush-Kuhn-Tucker conditions in each linear section. The algorithm for solving~\eqref{eq:proj_1inf_ball1_turlach} involves starting at $0$ and finding the ``knots'' where the $\rho_l$ change from one linear piece to another, until they encounter a interval containing $\tau$. 

 \subsection{Solution via linear programming}
 \label{sec:quattoni}
An equivalent linear program given by:

 \begin{align}
 	{\text{find}} & \quad \mathbf{\mu}, \theta \label{eq:proj_1inf_quattoni} \notag\\ 
\text{ s.t. } & \sum_n \mu_n = \tau \;,\notag \\
	& \sum_{m} (A_{nm}-\mu_n)_+ = \theta, \forall n \; \text{such that} \; \mu_n >0 \;,\notag \\
    & \sum_{m} A_{i,j} \leq \theta \;, \forall n \; \text{such that} \; \mu_n = 0 \;,\notag \\
    & \mu_n \geq 0, \forall n\;, \notag \\
    & \theta \geq 0 \;,
 \end{align}
was derived in \cite{quattoni2009efficient}.
The variables $\mu_i$ correspond to the $\ell_{\infty}$-norm values of each row of the optimum solution, and $\theta$ is related to a shrinkage parameter associated with the projection onto the $\ell_1$-ball. $\theta$ and $\mu$ are found via a search procedure over a piece-wise linear function, similar to that of Section~\ref{sec:turlach}.

\subsection{Solution via a dual splitting of the maximum function}
\label{sec:dualMax}
{The segmentation of an image into $K$ non-overlapping regions can be posed as a spatially regularized version 
of the K-means problem, i.e. piecewise constant Mumford-Shah problem. Via a lifting-based reformulation 
\cite[Section 2.1]{condat-2017-convex}, this problem can be expressed as
\begin{align}
 \min_{Z} \langle Z,\, W \rangle + \frac{\lambda}{2} \sum_{m=1}^M \text{TV}(\mathbf{z}_m) \; \mbox{ s.t. } \| Z \|_{{\infty,1}} \leq K,
 \label{eqn:condat}
\end{align}
where $\langle \cdot\, , \cdot \rangle$ represents the inner product,  TV$(\cdot)$ is any discrete form of 
Total Variation (TV) \cite{rudin-1992-nonlinear}, $\lambda$ is the regularization parameter, and $Z$ represents 
the set of $M$ groups, individually represented by $\mathbf{z}_m$, each of which as a binary assignment
vector, with elements in $\{0, 1\}~$\cite[Section 2.1]{condat-2017-convex}.

In order to avoid a direct solution of the projection onto the $\ell_{\infty,1}$ ball subproblem in (\ref{eqn:condat}), 
\cite{condat-2017-convex} made used of the infimal convolution (or epi-sum) representation of the maximum 
function \cite{bauschke-2011-convex}, \cite[Eq. (11)]{condat-2017-convex},
resulting in an splitting approach~\cite[Eq. (12)]{condat-2017-convex} tailored
to the particular properties of variable $Z$. While an adaptation of this method could, in
principle, be used to solve \eqref{eq:proj_1inf_ball1}, such an adaptation is not straightforward, and to the best of our knowledge,
has yet to be developed.
}

\subsection{Projection onto the $\ell_{p,1}$ ball by root search}
\label{sec:sra}

The general $\ell_{p,1}$-ball projection problem was solved by means of a root search technique in \cite{sra2011fast}. 
This approach relies on the fact that the proximal operator of the $\ell_{p,1}$ norm, defined as 
\begin{equation}
\mathrm{prox}_{\|\cdot\|_{p,1}}(B,\lambda)  := 
\argmin_X \frac{1}{2} \| X - B  \|_F^2 +
 \lambda \| X \|_{p,1} \;,
\label{eq:prox_inf_q}
\end{equation}
has a simpler solution than the projection onto the $\ell_{p,1}$-ball (\ref{eq:proj_1inf_ball1}),
%
since (\ref{eq:prox_inf_q}) can be computed by solving independent $\ell_p$-norm proximity subproblems~\cite{sra2011fast} of the form 
\begin{equation}
\argmin_{\mathbf{x}_m} \frac{1}{2} \| \mathbf{x}_m - \mathbf{b}_m  \|_F^2 + \lambda  \| \mathbf{x}_m \|_{p}\;.
\end{equation}

One of the contributions of \cite{sra2011fast} was to propose a method to take advantage of the separability of (\ref{eq:prox_inf_q}) in order to solve (\ref{eq:proj_1inf_ball1}). Let ${\cal{L}}(X, \theta)$ be the Lagrangian of (\ref{eq:proj_1inf_ball1}), i.e, 
\begin{equation}
	{\cal{L}}(X, \theta) = \frac{1}{2} \| X - B  \|_F^2 +  \theta \left( \| X \|_{\infty,1} - \tau \right)
\end{equation}
and let $\theta^*$ be the optimal dual variable. As long as $\tau>0$,~\eqref{eq:proj_1inf_ball1} satisfies Slater's conditions for strong duality~\cite{Boyd2004book}. Therefore, the primal optimal variable $X^*(\theta^*) = {\text{argmin}}_{X} \; \mathcal{L}(X,\theta^*)$ can be obtained by computing
\begin{equation}
      X^*(\theta^*) = \argmin_X \frac{1}{2}\| X - B  \|_F^2  +
      \theta^*( \| X \|_{p,1} -  \tau) \;. \label{eq:lagrangian}
\end{equation}
Defining
\begin{equation}
X(\theta) = \mathrm{prox}_{\|\cdot\|_{p,1}}(B,\theta)
\end{equation}
and
scalar function
\begin{equation}
      g(\theta) = \| X(\theta) \|_{p,1} -\tau \;,
      \label{eq:search_fun}
\end{equation}
it can be shown~\cite[Lemma 2]{sra2011fast}~\cite[Lemma 1]{sra-2012-fast} that there exist an interval $[0,\, \theta_{\text{max}}]$ over which
$g(\theta)$ is monotonically decreasing and differs in sign at the endpoints. Since $\theta^*$ coincides with the unique root of $g(\theta)$, it can be found by using a root finding method.
%

{A root finding method combining bisection, inverse quadratic interpolation, and the secant method was proposed in~\cite{sra2011fast}.} This root finding based solution for projections onto the $\ell_{p,1}$ ball was extended to the more general $\ell_{p,q}$ case in a follow-up article \cite{sra-2012-fast}.
Computational performance comparisons indicated \cite[Section 3.1]{sra-2012-fast} that this algorithm was both much more accurate and twice as fast as that of \cite{quattoni2009efficient} for the particular case of projections onto the $\ell_{\infty,1}$ ball.

\subsection{Solution using Steffensen's root-finding method}
\label{sec:steff}

We have previously presented a technique for solving the $\ell_{1,\infty}$-constraint problem via a Steffensen root-finding 
method~\cite{chau2018steffensen}. Steffensen's method is a quasi-Newton root finding algorithm~\cite{amat2016overview} that is useful when an analytical expression of the derivative is not available. The drawback is that it requires two function evaluations, and is usually more expensive than Newton's method.  Given the function $f(x)$, Steffensen's original iterations consist of the update
\begin{equation}
      x_{n+1} := x_n + \frac{x_n}{\delta F(x_n,y_n)} \;,
\end{equation}
where
\begin{align}
\delta F(x_n,y_n) &= \frac{f(y_n)-f(x_n)}{y_n - x_n} \\
y_n &= x_n + f(x_n) \;.
\end{align}

Steffensen's method tends to exhibit convergence problems if the initial $x_0$ is too far from the actual root. Therefore, \cite{chau2018steffensen} used the modified version proposed in~\cite{amat2016overview}:
\[
y_n = x_n + \alpha_n |f(x_n)| \;.
\]
Here $\alpha_n$ is an adaptive parameter that is recommended to take values that satisfy
\begin{equation}
    \mathrm{tol}_c \ll \frac{ \mathrm{tol}_u}{2|f(x_n)|} < |\alpha_n| < \frac{ \mathrm{tol}_u }{|f(x_n)|} \;,
    \label{eqn:st-tol}
\end{equation}
where $\mathrm{tol}_c$ is chosen in accordance 
with the numerical precision used in the implementation, and $\mathrm{tol}_u$ is a user-defined parameter.

To the best of our knowledge, the {fastest methods for solving the $\ell_{\infty,1}$-ball projection problem are~\cite{quattoni2009efficient},  \cite{sra2011fast} and~\cite{chau2018steffensen}. Accordingly, these three algorithms will be used as benchmark for all our comparisons.}

\section{Preliminaries}
\label{sec:preliminaries}
\subsection{Notation}

We will denote matrices with non-bold upper case font and vectors with bold lower case. Additionally, if $A \in \mathbb{R}^{M \times N}$ is a matrix, we will denote 
the $m^{\text{th}}$ row of $A$ by $\mathbf{a}_m$, and the $i^{\text{th}}$ element of the $m^{\text{th}}$ row of $A$ by $a_{im}$. 

\subsection{Projection onto the $\ell_1$-ball}
\label{sec:l1ball}

For our proposed solution to (\ref{eq:proj_1inf_ball1}), we will need to solve the closely related problem of projection onto the $\ell_1$-ball, which is defined as
 \begin{equation}
  \mathrm{proj}_{\|\cdot\|_1}(\mathbf{u},\tau) = \underset{X} \argmin \; \frac{1}{2} \| \mathbf{x} -  \mathbf{u} \|_2^2 \;\;
  \mbox{ s.t. } \| \mathbf{x} \|_1 \leq \tau \;,
\label{eq:ball1}
\end{equation}
where $\mathbf{x}, \mathbf{u} \in \mathbb{R}^N$. 
%
The solution to this equation is given by \cite{duchi-2008-efficient, liu-2009-efficient, sra-2010-generalized,
sra-2012-fast, condat2016fast}
%
%
\vspace{-2mm}
\begin{equation}
  \mathbf{x}^*=
\begin{cases}
  \mathbf{u} & \text{if}\ \| \mathbf{u} \|_1 < \tau  \\
  \shrink(\mathbf{x}, \lambda(\tau)) & \text{if}\ \| \mathbf{u} \|_1 \geq \tau \;,
\end{cases}
\label{eq:l1ball_cases}
\end{equation}
where 
\begin{equation}
	\shrink(\mathbf{x}, \lambda(\tau)) = \text{sign}(\mathbf{x})\odot \max (|\mathbf{x}|-\lambda(\tau),0) \;,
\end{equation}
$\lambda(\tau)$ is a shrinkage parameter that depend on $\tau$, and $\odot$ is the element-wise (Hadamard) vector product.

To solve problem  (\ref{eq:ball1}), we must find $\lambda^*$ such that $f(\lambda^*) = 0$, where
\begin{equation}
    f(\lambda) = \sum_n \max(|u_n|-\lambda,0)-\tau \;.
    \label{eq:pl1-func}
\end{equation}

Clearly, if $|u_n| > \lambda^*$ then this element contributes to the sum defined in (\ref{eq:pl1-func}); thus, 
by sorting $| \mathbf{u} |$ in decreasing order, a simple search will lead to the solution:
i.e. let $\mathbf{v} = \mbox{\texttt{sort}}( | \mathbf{u} | )$ and define 
\begin{equation}
	L = \max\left\{ l \, \Bigg| \, l^{-1} \left(\sum_{n=0}^l v_n - \tau\right) < v_l \right\}\;, 
\end{equation}
then
\[
\lambda^* = L^{-1} \left(\sum_{n=0}^L v_k - \tau \right)  \;.
\]
This solution was originally described 
in \cite{held-1974-validation}, with several improvements also reported in 
\cite{berg-2009-probing,kiwiel-2008-breakpoint,duchi-2008-efficient}.

The solution of problem  \eqref{eq:ball1} can be cast as a root-finding problem, as the function defined in (\ref{eq:pl1-func}) satisfies the Fourier conditions in the interval $[0, u_{\mbox{max}}]$ and thus has a single root in this interval~\cite{liu-2009-efficient,gong-2011-efficient}.

The application of the Newton root-finding method~\cite[Section 11.1]{NoceWrig06}: i.e.
\[
\lambda_{k+1} = \lambda_k - \frac{f(\lambda_k)}{f'(\lambda_k)}
\]
leads to the so-called Michelot algorithm \cite{michelot-1986-finite,kiwiel-2008-variable,liu-2009-efficient, gong-2011-efficient}. More recently, \cite{condat2016fast,rodriguez2017parallel,rodriguez-2018-accelerated} have
proposed further improvements to \cite{michelot-1986-finite}.


Since it will be helpful in the 
derivation of our proposed algorithm (see Section \ref{sec:proposed}), we show how the Michelot algorithm can be derived from a Newton's root-search method applied to \eqref{eq:pl1-func}. This equation can be rewritten as~\cite[eq. (32)]{gong-2011-efficient},\cite[eq. (7)]{rodriguez2017parallel}:
\begin{equation}
    f(\lambda) = \mathbf{z^Tu} - \lambda\mathbf{z^Tz} - \tau\;, \label{eq:l1search}
\end{equation}
where
\begin{equation}
    \mathbf{z} = \text{sign}(\mathbf{u}) \odot I_{|\mathbf{u}|<\lambda}\;,
\end{equation}
where $I(\cdot)$ is the indicator function defined for a set $A$ as \cite[Chapter 2]{folland2013real}:
\begin{equation}
	I_A(x) = \begin{cases}
		1 &  x \in A \\
        0 &  x \notin A.
	\end{cases}
\end{equation}

Applying the Newton method to (\ref{eq:l1search}) by temporarily disregarding the dependence of $\mathbf{z}$ on $\lambda$, leads to the Michelot algorithm iterations:
\begin{equation}
 \lambda_{k+1} =  
 \lambda_k - 
 \frac{ \mathbf{z}_k^T\mathbf{u} - \lambda_k\mathbf{z}_k^T\mathbf{z}_k - \tau }{-\mathbf{z}_k^T\mathbf{z}_k}  
  \;=\;  \frac{ \mathbf{z}_k^T\mathbf{u} - \tau}{ \mathbf{z}_k^T\mathbf{z}_k }\;,
 \label{eq:mich-new1}
\end{equation}
where  $\mathbf{z}_k^T\mathbf{u}$ is the $\ell_1$ norm of the subset of elements of $mathbf{u}$ 
for which the corresponding absolute value is greater than $\lambda_k$,
and $ \mathbf{z}_k^T\mathbf{z}_k $ is the number of elements of this subset (equal to the 
number of non-zero elements in $ \mathbf{z}_k$).
Furthermore it can also be shown that
\begin{equation}
 0 < \lambda_{k} \leq \lambda_{k+1} \leq \lambda^* \; \forall k\;,
 \label{eqn:lambda}
\end{equation}
where $\lambda^*$ is the parameter that solves (\ref{eq:ball1}).
Thus  at each iteration we can discard or prune all the elements in $\mathbf{u}$
such that $| u_n | \leq \lambda_k$. Following the guidelines given in \cite[Section 3.2.2]{rodriguez2017parallel},
a careful implementation of such a pruning strategy can lead to a very efficient computational performance of the Michelot algorithm, as empirically shown in \cite[Tables I and II]{rodriguez-2018-accelerated}.

\subsection{Dual norm}
\label{sec:dual}
In this section, we summarize additional theoretical results from \cite{sra2011fast} that will be useful in the analysis and derivation of the proposed algorithm.
\begin{definition} \cite{Boyd2004book} If $\| \cdot \|$ is a norm on $\mathbb{R}^m$, then the associated dual norm, $\| \cdot \|_*$, is defined as
\begin{equation}
	\| \mathbf{z} \|_* \triangleq \sup \left\lbrace \mathbf{z}^T \mathbf{x} \; | \; \| \mathbf{x}\| \leq 1 \right\rbrace \;.
\end{equation}
\end{definition}
%

\begin{lemma}[see {\cite[Lemma 1]{sra2011fast}}]
  Let $q \geq 1$ and let $p^*$ be its conjugate exponent satisfying $\frac{1}{p} + \frac{1}{p^*} = 1$. Then, the norm $\| \cdot \|_{p^*,\infty}$ is dual to $\| \cdot \|_{p,1}$.
\label{lemmadual}
\end{lemma}

It can be shown via Moreau's decomposition \cite{parikh2014proximal}, \cite{sra2011fast} that the dual problem of
\begin{equation}
 \mathrm{prox}_{\|\cdot\|_{p^*,\infty}}(B, \lambda) :=
 \argmin_X \; \frac{1}{2} \| X - B  \|_F^2 +
 \lambda \cdot \| X \|_{p^*,\infty} \;,
 \label{eqn:prox-infty1}
\end{equation}
where
\[
\| X \|_{p^*,\infty} = \max \{ \| \mathbf{x}_k \|_{p^*} \} \;,
\]
 is the projection onto the $\ell_{p,1}$ ball with radius $\lambda$, i.e., 
 \begin{align}
\mathrm{proj}_{\|\cdot\|_{p,1}}(B,\lambda) := & \;
\argmin_X
\frac{1}{2} \| X - B  \|_F^2 \label{eq:proj_p1ball} \\ 
& \text{ s.t. } \| X \|_{p,1} \leq \lambda\;. \notag
\end{align}

\section{Proposed method}
\label{sec:proposed}

\subsection{Leveraging $\mathrm{prox}_{\|\cdot\|_{1,\infty}}(\cdot)$, the dual of
$\mathrm{proj}_{\|\cdot\|_{\infty,1}}(\cdot)$}
\label{sec:sol_prox}
%

As described in Section \ref{sec:sra}, the approach of~\cite{sra2011fast} involved solving~(\ref{eq:proj_p1ball}) via a root finding method applied to (\ref{eq:proj_1inf_ball1}). Here we consider an alternative reinterpretation that allows us to derive several improvements in the proposed algorithm, as explained in the following section.
By the results of Section~\ref{sec:dual}, the proximal operator of $\ell_{1,\infty}$ is the dual of the projection on the $\ell_{\infty,1}$ ball and vice-versa, then $X^*=\mathrm{proj}_{\|\cdot\|_{\infty,1}}(B, \tau)$, can be written as $X^* = B-A^*$, where 
%
\begin{align}
    A^* &= \mathrm{prox}_{\|\cdot\|_{1,\infty}}(B,\tau) \nonumber \\
        &= \argmin_A   \frac{1}{2} \| A - B  \|_F^2 +
    \tau \| A \|_{1,\infty} \label{eq:prox_inf1} \;.
 \end{align}
%
%
%
Now, if $A^*$ is known, we can define 
\[
\gamma^* = \| A^* \|_{1,\infty} = \max_m \{ \| \mathbf{a}^*_m \|_1 \} \;,
\]
and thus, after simple algebraic manipulation,~\eqref{eq:prox_inf1} can be written as 
%
\begin{equation}
\argmin_{\{\mb{a}_m\}}
 \; \frac{1}{2} \sum_m \| \mathbf{a}_m- \mathbf{b}_m  \|_2^2 
 \quad \text{ s.t.} \quad \| \mathbf{a}_m \|_1 \leq \gamma^*,\, \forall m \;. \label{eq:proj_1inf_ball} 
\end{equation}



Clearly,~\eqref{eq:proj_1inf_ball} is separable in $\mathbf{a}_m$, with the individual problems corresponding to a projection on the $\ell_1$-ball (see Section \ref{sec:l1ball}).
Accordingly, if we devise a method for obtaining the optimal $\gamma^*$ value, then the solution to~\eqref{eq:prox_inf1}, and therefore to~\eqref{eq:proj_1inf_ball1}, can be easily calculated. The $\gamma^*$ value can be found by a root finding method, as described in the following section.

\subsection{Search function and solution by Newton's method}
\label{sec:search_fun}
As originally proposed in \cite{sra2011fast}, we use $\mathrm{prox}_{\|\cdot\|_{\infty,1}}(\cdot)$ to solve $\mathrm{proj}_{\|\cdot\|_{\infty,1}}(\cdot)$ in (\ref{eq:search_fun}). Thus, we replace $X$ by $B-A$ in~\eqref{eq:search_fun} and after simple algebraic manipulations, we obtain
\begin{eqnarray}
 f(\gamma) &=& \sum_{m=1}^M {||\mathbf{b}_m - \mathbf{a}_m(\gamma)||_{\infty} - \tau} \;, \label{eq:search_prop}
\end{eqnarray}
defined for $\gamma \geq 0$. Furthermore, since (\ref{eq:search_prop}) is equivalent to (\ref{eq:search_fun}),
it also satisfies the Fourier conditions, and thus
it has a unique root at $\gamma^*$. For a given $\gamma$, $\mathbf{a}_m$ is computed using the approach
described in~\eqref{eq:proj_1inf_ball}.
As each $\mathbf{a}_m(\gamma)$ corresponds to a projection onto the $\ell_1$-ball,
we apply~\eqref{eq:l1ball_cases} and obtain

\begin{small}
\begin{equation}
  \mathbf{a}_m(\gamma) =
\begin{cases}
  \mathbf{b}_m & \text{if}\ \| \mathbf{b}_m \|_1 < \gamma  \\
  \text{sign}(\mathbf{b}_m)\odot \max (|\mathbf{b}_m|-\lambda_m(\gamma),0) & \text{if}\ \| \mathbf{b}_m \|_1 \geq \gamma\;.
\end{cases}
\label{eq:l1ball_cases2}
 \end{equation}
\end{small}

%

%

By substituting~\eqref{eq:l1ball_cases2} into~\eqref{eq:search_prop}, we note that only the terms corresponding to the $\| \mathbf{b}_m \|_1 \geq \tau$ contribute to the sum. Accordingly, at each evaluation of the search function, we can prune the rows of $B$ that do not fulfill this condition, and only perform the projections specified in~\eqref{eq:l1ball_cases2} on the
remaining rows. Our numerical experiments show that this pruning strategy can reduce the computational time by half or more. 
Based on this remark, we can rewrite the search function as:

\begin{small}
\begin{equation}
  f(\gamma) = \sum_{m \in \mathcal{M}} \hspace{-1mm}||\mathbf{b}_m -  
   \text{sign}(\mathbf{b}_m)\odot \max (|\mathbf{b}_m|-\lambda_m(\gamma),0) ||_{\infty} - \tau \;, \label{eq:search_reduced}
\end{equation} 
\end{small}
%
%

{\noindent}where $\mathcal{M}$ denotes the set of indexes $m$ where $\| b_m\|_1 \geq \gamma$. 
We can reduce this expression further by noting that $b_m$ can be rewritten in the 
form $\text{sign}(\mathbf{b}_m)\odot |\mathbf{b}_m|$ and factorizing:

 \begin{small}
  \begin{align}
  f(\gamma)  = &   \sum_{m \in \mathcal{M}} ||\text{sign}(\mathbf{b}_m)\odot 
   (|\mathbf{b}_m| - \max (|\mathbf{b}_m|-\lambda_m(\gamma),0))  ||_{\infty} - \lambda \\
  f(\gamma)  = &\sum_{m \in \mathcal{M}} {| \; ||\mathbf{b}_m| - \max (|\mathbf{b}_m|-\lambda_m(\gamma),0)  ||_{\infty} - \tau} \;. \label{eq:summax}
\end{align}
\end{small}

 Now, we turn to the analysis of 
 \begin{equation}
 \beta(\mathbf{b}_m) = || \; |\mathbf{b}_m| - \max (|\mathbf{b}_m|-\lambda_m(\gamma),0)  ||_{\infty}\;.
 \end{equation}

   We will denote by $b_{im}$ the $i^\text{th}$ component of the vector $\mathbf{b}_m$. As all the components of this vector are positive, we can write 
  \begin{equation}
  	\beta(\mathbf{b}_m) = \max_i (|b_{im}| - \max (|b_{im}|-\lambda_m(\gamma),0)) \;.
  \end{equation}
  For each component, we have:
   \begin{equation}
  \small
   |b_{im}| - \max (|b_{im}|-\lambda_m(\gamma),0)= 
\begin{cases}
  \lambda_m(\gamma) & \text{if}\ | b_{im} | > \lambda_m(\gamma)  \\
  | b_{im} | & \text{if}\ | b_{im} | \leq \lambda_m(\gamma)
\end{cases} \;.
\label{eq:l1ball_cases3}
 \end{equation} 

Now, we assert that there exists at least one element $b_{jm}$ such that $| b_{jm} | > \lambda_m(\gamma)$. To prove this, suppose that $| b_{im} | \leq \lambda_m(\gamma)$ for all $i$. Substituting this assumption into \eqref{eq:l1ball_cases}, we would have that $\mathbf{a}_m^*$ is zero. As we are only considering terms corresponding to $\| \mathbf{b}_m \|_1 \geq \gamma$ we arrive to a contradiction.

Then, as there exist $b_{jm}$ such that $| b_{jm} | > \lambda_m(\gamma)$

\[
|b_{jm}| - \max (|b_{jm}|-\lambda_m(\gamma),0)=\lambda_m(\gamma) \;.
\]
All other elements are in turn less or equal to $\lambda_m(\gamma)$. From this, we conclude that  $\beta(\mathbf{b}_m) = \lambda_m(\gamma)$. Thus, we rewrite~\eqref{eq:summax} as
  \begin{equation}
  f(\gamma) = \sum_{m \in \mathcal{M}} {\lambda_m(\gamma) - \tau} \label{eq:search_minexp}\;.
 \end{equation}
 
 As outlined in equation~\eqref{eq:mich-new1} (note that here the sub-indexes have a different interpretation), 
 $\lambda_m(\gamma)$ can be expressed as
 \begin{equation}
	\lambda_m(\gamma) = \frac{z_k^T \mathbf{b}_m - \gamma}{z_k^T z_k}\;,
 \end{equation}
 where
 \[
 z_m = \text{sign}(\mathbf{b}_m) \odot I_{|\mathbf{b}_m|<\lambda_m(\gamma)} \;,
\]
 so that \eqref{eq:search_minexp} becomes
\begin{equation}
  f(\gamma) = \sum_{m \in \mathcal{M}} {\frac{z_m^T \mathbf{b}_m - \lambda}{z_m^T z_m} }- \tau\;.
 \end{equation}
 
Both $z_m$ and $\mathcal{M}$ depend on $\gamma$, however, similar to the derivation for the 
Michelot algorithm for $\ell_1$-ball projection~\cite{cominetti2014newton} presented in Section~\ref{sec:l1ball}, we temporarily disregard these dependencies and approximate the derivative of $f$ as
 \begin{equation}
   \frac{\partial f(\gamma)}{\partial \gamma} \approx -\sum_{m \in \mathcal{M}}{ \frac{1}{z_m^T z_m} }\;. \label{eq:derivative}
 \end{equation}

 Thus, the updates of the root-finding procedure can be performed in a Newton-like fashion by setting:
 
 \begin{equation}
        \gamma_{n+1} := \gamma_n + \frac{f(\gamma)}{\sum_{m \in \mathcal{M}}{ \frac{1}{z_m^T z_m}}}   \;.
 \end{equation}

  Our numerical experiments suggest that, if we update $z_k$ at each iteration, this approximation is good enough for use in a Newton root search method. Similarly to the re-derivation of the Michelot algorithm presented in ~\cite{cominetti2014newton}, our method can be understood as a quasi-Newton method in the broad-sense of the term. However, we note that it cannot be readily derived from the application of classical quasi-Newton schemes for root-finding, such as Broyden's or Brent's method~\cite{Chapra2006book}.

\subsection{Initial Point}
\label{sec:ini_point}

From here on, we suppose that
\begin{equation}
\| B \|_{\infty,1} = \sum_{m=1}^M {\| \mathbf{b}_m \|_{\infty}} > \tau \;.
\end{equation}
If $\| B \|_{\infty,1} \leq \tau$ in~\eqref{eq:proj_1inf_ball1} then the optimal solution is trivial, $X^*=B$. We try to find a point $\gamma_0$ such that $f(\gamma_0)>0$ in~\eqref{eq:search_prop}. Then, as $f(\cdot)$ satisfies the Fourier conditions and is therefore non-increasing in the $[0,\gamma^*]$ interval, we can conclude that $0 \leq \gamma_0 \leq \gamma^*$.

We start by assuming that the $\ell_1$-norm of the $j^{\text{th}}$ row of the solution $\mathbf{A^*}$ coincides with
$\| \mathbf{A^*} \|_{1,\infty}$, i.e., $\max_m \{ \| \mathbf{a}_m \|_1 \} = \| \mathbf{a}_j \|_1$. Then,
via (\ref{eq:prox_inf1}),
we can find $\mathbf{a}_j$ as:
 \begin{equation}
       \mathbf{a}_j =  \argmin_{\mb{a}} \frac{1}{2} \| \mathbf{a} - \mathbf{b}_j \|^2_2
       + \tau \| \mathbf{a} \|_1 = \text{shrink}(\mathbf{b}_j,\tau) \;. \label{eq:shrink}
 \end{equation}
 
We define $\gamma_0 = \| \text{shrink}(\mathbf{b}_j,\tau) \|_1$
and proceed to show that $f(\gamma_0)>0$. 
 Separating the sum and using the definition of the shrinkage operator in~\eqref{eq:search_prop}, we can write:
 
   \begin{small}
 \begin{align}
  f(\gamma_0)  = & \sum_{m=1}^M {||\mathbf{b}_m - \mathbf{a}_m(\gamma)||_{\infty} - \tau} \label{eq:search_shrink} \\
  f(\gamma_0)  = & \| b_j-\text{sign}(b_j)\odot (|b_j| - \max (|b_j|-\tau,0))  \|_{\infty} \\ \notag
    + & \sum_{m \neq j} {||\mathbf{b}_m - \text{proj}_{\| \cdot \|_1}(b_m,\gamma_0)||_{\infty} - \tau} \;. \label{eq:search_div}
\end{align}
 \end{small}
 
 We can reduce this expression further by noting that $b_j$ can be rewritten in the form $\text{sign}(b_k)\odot |b_k|$ and factorize:
\begin{equation}
  f(\gamma_0)  =  \tau_j(\mathbf{b}_j) 
  +  \sum_{m \neq j} {||\mathbf{b}_m - \text{proj}_{\| \cdot \|_1}(\mathbf{b}_m,\gamma_0)||_{\infty} - \tau} \;, \label{eq:search_div_2}
\end{equation}
where 
\begin{equation}
	\tau_j(\mathbf{b}_j) = \| \; |\mathbf{b}_j| - \max (|\mathbf{b}_j|-\tau,0) \|_{\infty} \;.
\end{equation}
  
  We now turn to the analysis of $\tau_j(\mathbf{b}_j)$.  As all the components involved in $\tau_j(\mathbf{b}_j)$ are positive, we can write the norm as the maximum of all the components of the vector. For each component, we have:
\begin{equation}
  |{b}_j^{(i)}| - \max (|{b}_j^{(i)}|-\tau,0) =
\begin{cases}
  \tau & \text{if}\ | b_{ij} | > \tau  \\
  | b_{ij} | & \text{if}\ | b_{ij} | \leq \tau \;.
\end{cases}
\label{eq:cases_shrink}
\end{equation}

If we assume that $\| \mathbf{b}_j \|_{\infty}>\tau$, as this is a simple algorithmic check included in our method (See Section~\ref{sec:algoritmo}), then at least one of the components of $\mathbf{b}_j$ must fulfill the first condition in~\eqref{eq:cases_shrink}. For at least this component, we have 
\begin{equation}
	|b_{ij}| - \max (|b_{ij}|-\tau,0) = \tau\;,
\end{equation}
and all the other components are less than or equal to $\tau$. Accordingly,
\[
\| \; |\mathbf{b}_j| - \max (|\mathbf{b}_j|-\tau,0) \|_{\infty} = \tau \;.
\]
Replacing this in~\eqref{eq:search_div_2}, we obtain
 \begin{equation}
  f(\gamma_0) = \sum_{m \neq j} {||\mathbf{b}_m - \text{proj}_{\| \cdot \|_1}(\mathbf{b}_m,\gamma_0)||_{\infty}} > 0 \;.
\end{equation}

Thus, $f(\gamma_0)>0$ and, instead of starting the root search from 0, we can start from $\gamma_0$, which is a better initial guess of $\gamma^*$.  \textit{A priori}, we do not know which $j$ is closer to the real maximum.
In order to obtain the initial point $\gamma_0$,
we solve~\eqref{eq:shrink} for every row and then take among these solutions the one with the maximum $\ell_1$-norm.

\subsection{Proposed method}
\label{sec:algoritmo}

\begin{figure*}
  \centering
       \includegraphics[width=\textwidth]{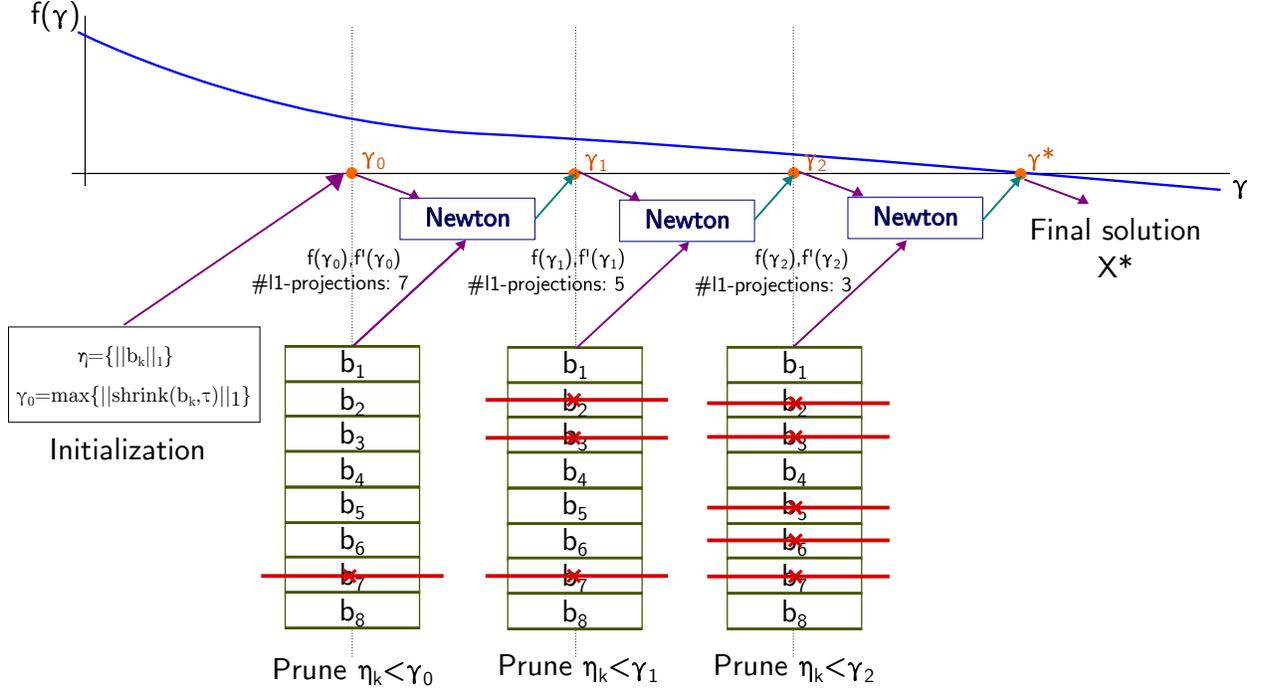}
      \vspace{-1mm}\caption{Diagram illustrating the proposed method for an specific example. 
      The initialization block provides the initial point $\gamma_0$ for the root search procedure (see \ref{sec:initial}). At each iteration of $\gamma_n$ the matrix rows are pruned based on their $\ell_1$-norm as shown in the bottom part of the graphic. The updates of $\gamma_n$ are done via Newton's method and they converge to the root of the function as shown in the upper part of the graphic (see \ref{sec:search_fun}) }
      \label{fig:block_diag}
\end{figure*}

The full proposed method is presented in Algorithm~\ref{al:proposed} and depicted for an specific example in Figure~\ref{fig:block_diag}.
In line \ref{algline:menor}, $B$ is checked to see if it is already in the $\ell_{\infty,1}$-ball.
Lines \ref{algline:ini1} and \ref{algline:ini2} compute the initial guess of the solution as outlined in Section~\ref{sec:ini_point} and the initialization Block of Figure~\ref{fig:block_diag}.
Note that if $\| B \|_{\infty,\infty} = \max_{i,j} b_{i,j} < \tau$, lines \ref{algline:ini1} and \ref{algline:ini2} do not need to be evaluated
and $\gamma$ is assigned an initial value of 0. Likewise, the shrinkage operation is performed
solely for the rows whose $\ell_{\infty}$-norm is greater than $\tau$.

Line \ref{algline:prune} corresponds to the pruning step described in Section~\ref{sec:search_fun}, where all rows of $B$ with $\ell_1$ norm less than the current $\gamma$ are discarded. This is also illustrated in the bottom part of the graphic, where the $\ell_1$-projections are not performed on the discarded rows. Lines \ref{algline:der} and \ref{algline:update} perform the approximation of the derivative and the Newton updates (see Section~\ref{sec:search_fun}, Blue Newton blocks in the Figure~\ref{fig:block_diag}). Line \ref{algline:final} obtains the final $\mathbf{A^*}$ with the last updated $\gamma$ value and Line \ref{algline:moreau} returns $B-\mathbf{A^*}$ due to Moreau's decomposition (see Section~\ref{sec:sol_prox}).

\begin{algorithm}
\caption{Proposed method via root-finding}
      \label{al:proposed}
      \hspace*{\algorithmicindent} 
       \begin{algorithmic}[1]
       \IF{$\| B \|_{\infty,1} \leq \tau$}
       		\STATE{return $B$} \label{algline:menor} 
       \ENDIF
 \IF{$\| B \|_{\infty,\infty}>\tau$} 
 		\STATE{Compute $\alpha_k = \| \mathrm{shrink} (b_k,\tau) \|_1$ for each row of $B$.}  \label{algline:ini1}
        \STATE{Define $\gamma = \max_k (\alpha_k)$ } \label{algline:ini2}
       \ELSE
\STATE{$\gamma = 0$.}
	\ENDIF
    \FOR{$k=1:\mathrm{maxIter}$}
    	\STATE{ Prune the rows of $B$ that have $\ell_1$-norm less than $\gamma$} \label{algline:prune}
        \STATE{Obtain $f(\gamma)$ as defined in~\eqref{eq:search_prop} and $\frac{\partial f(\gamma)}{\partial \gamma}$ as defined in~\eqref{eq:derivative}. \label{algline:der}}
        \IF{$|f(\gamma)|< \mathrm{tolerance}$}
           \STATE{\textbf{break}}
        \ENDIF
        \STATE{Update $\gamma$ using Newton method. \label{algline:update}}
    \ENDFOR
\STATE{Solve for $A$ in~\eqref{eq:proj_1inf_ball} with the obtained $\gamma$.}     \label{algline:final}
\STATE{Return $B-A$} \label{algline:moreau}
     
\end{algorithmic}

\end{algorithm}

%
%

\begin{algorithm}
\caption{Projected gradient descent}
      \label{al:pgd}
      \hspace*{\algorithmicindent} \textbf{Input:} matrix $B$, projectionOperator, maxIter, tolerance, $\alpha$, $W_0$ \\
      \hspace*{\algorithmicindent} \textbf{Initialization:} $W := W_0$
       \begin{algorithmic}[1]
       \FOR{$k=1:maxIter$}
       	\STATE{Compute $A$ with columns $A^{(i)}:= W_{(k-1)}^{(i)} - \alpha X^T(XW^{(i)}-Y^{(i})$}
        \STATE{$W_{(k)}= \text{projectionOperator}(A)$}
        \IF{$\| W_{(k)}-W_{(k-1)}\|_F<$ tolerance}
        	\STATE{\textbf{break}}
        \ENDIF
       \ENDFOR
       \end{algorithmic}
  \end{algorithm}

\section{Multi-task LASSO}
\label{sec:mtl}

The Multi-task LASSO (MTL) problem will be used for testing our method on an application involving real data. 
Let $\text{vec}(\cdot)$ and $\text{vec}^{-1}(\cdot)$ be the vectorization operator and its inverse. Given $K$ tasks, 
each of length $N$ ordered in a $NK \times 1$ vector $\mathbf{b}$, and a coefficient matrix 
$W \in \mathbb{R}^{N \times M}$ we want to find the matrix $X \in \mathbb{R}^{M \times K}$ of features that solves
\begin{align}
\argmin_X & \; \frac{1}{2}  \| \mathbf{b}- P \text{vec}(X)  \|_2^2	  \notag \\ 
   \text{ subject to}  & \quad \| X \|_{1,\infty} \leq \tau \;,  \label{eq:mtl}
 \end{align}
where $P = I_N \otimes W$. The solution of this problem will tend to have few non-zero rows, i.e. selected features. It can be solved by means of projected gradient descent (PGD)~\cite[Chapter 3]{Bertsekas99book}, which is shown in Algorithm~\ref{al:pgd}. The method consists of alternating unconstrained gradient descent steps and projections into the $\| X \|_{1,\infty} \leq \tau$ ball. \eqref{eq:mtl} is a convex problem and thus it can be solved using other, more general, optimization methods such as interior point methods~\cite{Boyd2004book} or methods based on the augmented Lagrangian function~\cite{xu2017firstorder}.

In this case, projection operator will be a routine that solves~\eqref{eq:proj_1inf_ball1} by means of our proposed method or via one of the methods in the literature.


\section{Results}
\label{sec:rslts}

All tests presented below were computed using single-threaded Matlab or {C-Mex} code running on an 
Intel i7-4770K CPU (8 cores, {2.00} GHz, 32GB RAM). In our simulations with 
synthetic data (Sections \ref{sec:initial} and \ref{sec:simulations}), matrix $B$ was generated using a 
uniform distribution $[-0.5, 0.5]$, and $\tau$, the constraint used in (\ref{eqn:infty-1}), 
was taken such that $\tau = \alpha \| B \|_{\infty,1}$, where $\alpha$ is a small constant. 
Specific sizes of $B$ and values of $\alpha$ are mentioned below. 
Our Matlab {and C} code~\cite{chau-l1Infty-codes} can be used to reproduce our experimental results.

\subsection{Impact of initial point}
\label{sec:initial}

\begin{figure}[h]
  \centering
       \includegraphics[width=0.85\textwidth]{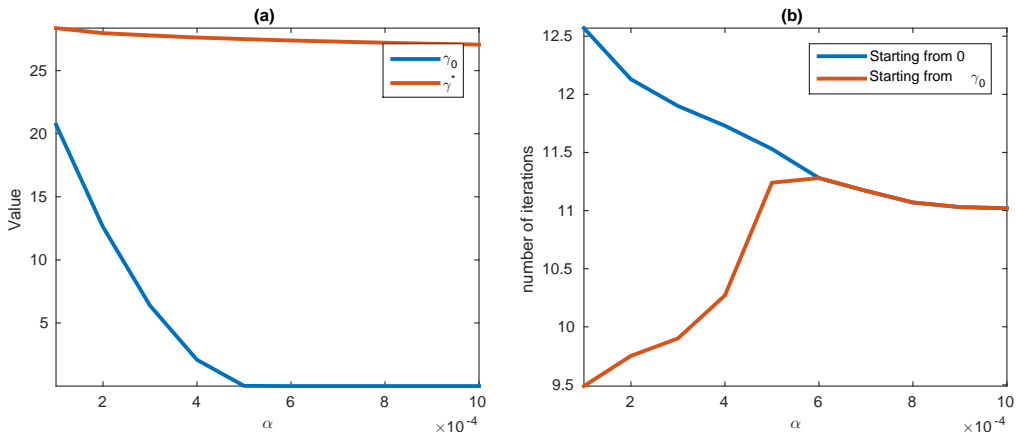}
      \vspace{-5mm}\caption{(a) Value of $\gamma_0$ (blue) and $\gamma^*$ (red) versus $\alpha$ (b) Number of
      iterations for proposed method to arrive at $\gamma^*$ starting from 0 (blue) and $\gamma_0$ (red)
      for different $\alpha$ values. See Section \ref{sec:initial}.}
      \label{fig:pinicial}
\end{figure}

In order to study the impact of the initial point $\gamma_0$ on the performance of the algorithm, 
we constructed 100 different realizations of a 2000 $\times$ 100 $B$ matrix, 
considering\footnote{These sizes and sparsity values are typical for known applications of (\ref{eqn:infty-1}) \cite{quattoni2009efficient, liu2009blockwise, sra2011fast}. Results for larger values of 
$\alpha$ can be obtained with our source code~\cite{chau-l1Infty-codes}.} $\alpha \in [10^{-4},\, 10^{-3}]$. 
For each value of $\tau$, the values for the initial point $\gamma_0$ and the optimal value $\gamma^*$ 
were averaged across the 100 realizations. These average values for each $\tau$ are shown in Figure~\ref{fig:pinicial}(a).  
It is observed that, at low $\alpha$ values, $\gamma_0$ is very close to the optimal value, but it goes rapidly to zero 
as $\tau$ increases. On the other hand, Figure~\ref{fig:pinicial}(b) 
shows a comparison of the average number of iterations that the proposed method needs 
for arriving to the optimal value starting from either zero or $\gamma_0$.  
The average number of iterations (across $100$ realizations)
and computational time for the different $\tau$ values, for the proposed method without 
pruning, along with the improvements provided by using $\gamma_0$, are listed in Table~\ref{tab:pinicial}.

\begin{table}
\addtolength{\tabcolsep}{-5pt}
\centering
\caption{ Computational results comparing the effect of the initial point $\gamma_0$.
  The percent change from the zero-start case is shown in parenthesis for the 
  $\gamma_0$ case. \textit{Num. Iter.} represents the average number of iterations across $100$ realizations. 
  See Section \ref{sec:initial}.}
\label{tab:pinicial}
\begin{tabular}{c|l|l|ll|ll}
 & \multicolumn{2}{c|}{\textbf{Starting at zero}} & \multicolumn{4}{c}{\textbf{Starting at $\gamma_0$}} \\ \hline
$\mathbf{\alpha\times 10^{-3}}$ / \textbf{sparsity(\%)} & \textbf{num iter} & \textbf{time(s)} & \multicolumn{2}{c|}{\textbf{num iter}} & \multicolumn{2}{c}{\textbf{time(s)}} \\ \hline
0.1 / 1.02 & 12.6 & 0.5 & 9.4 &(-25.0\%) & 0.18 &(-65.0\%) \\
0.2 / 1.92 & 12.2 & 0.5 & 9.7 &(-19.9\%) & 0.24 &(-52.4\%) \\
0.3 / 2.68 & 11.9 & 0.5 & 10.0 &(-16.7\%) & 0.28 &(-44.5\%) \\
0.4 / 3.40 & 11.7 & 0.5 & 10.2 &(-12.8\%) & 0.34 &(-32.2\%) \\
0.5 / 4.17 & 11.5 & 0.5 & 11.3 &(-2.0\%) & 0.48 &(-5.0\%) \\
0.6 / 4.94 & 11.3 & 0.5 & 11.3 &(0.0\%) & 0.52 &(0.3\%) \\
0.7 / 5.53 & 11.1 & 0.5 & 11.1 &(0.0\%) & 0.52 &(0.7\%) \\
0.8 / 6.28 & 11.1 & 0.5 & 11.1 &(0.0\%) & 0.52 &(0.7\%) \\
0.9 / 6.89 & 11.0 & 0.5 & 11.0 &(0.0\%) & 0.52 &(0.6\%) \\
1.0 / 7.53 & 11.0 & 0.5 & 11.0 &(0.0\%) & 0.52 &(0.6\%)
\end{tabular}
\end{table}

\subsection{Simulations}
\label{sec:simulations}

\begin{table}[]
\centering
\caption{Results for simulations with matrices of different size and the three tested methods 
in pure Matlab code. Error (Err.), number of iterations (N.I.) and running times are shown for 
each of them. Speedup with respect to GRF is shown for SRF and Proposed. Furthermore, we point out 
that the error of SRF and Proposed is three to five orders of magnitude lower than that of GRF.} 
\label{tab:comps}
\resizebox{\textwidth}{!}{{
\begin{tabular}{c|r|r|r|r||r|r|r|r||r|r|r|r}
 \textbf{} &  & \multicolumn{3}{c||}{\textbf{GRF}~\cite{sra2011fast} } & \multicolumn{4}{c||}{\textbf{SRF}~\cite{chau2018steffensen}} & \multicolumn{4}{c}{\textbf{Proposed}} \\ \hline
\textbf{Matrix Size} & \textbf{$\alpha$ / sp} & \textbf{Err.} & \textbf{N.I.} & \textbf{Time(s)} & \textbf{Err.} & \textbf{N.I.} & \textbf{Time(s)} & \textbf{Speedup} & \textbf{Err.} & \textbf{N.I.} & \textbf{Time(s)} & \textbf{Speedup} \\ \hline \hline
 & 0.0001/ 1.02 & 3.4e-11 & 9.6 & 9.2 & 2.6e-12 & 9.4 & 0.5 & \textbf{17.19} & 2.0e-16 & 9.4 & 0.3 & \textbf{33.54} \\ \cline{2-13} 
\textbf{2 000 x 100} & 0.0005/ 4.13 & 1.6e-10 & 13.2 & 9.5 & 1.4e-12 & 11.3 & 1.7 & \textbf{5.52} & 7.4e-16 & 11.3 & 0.7 & \textbf{13.02} \\ \cline{2-13} 
\textbf{} & 0.001/ 7.51 & 5.1e-10 & 14.5 & 9.7 & 1.3e-12 & 11.0 & 1.7 & \textbf{5.75} & 1.5e-15 & 11.0 & 0.8 & \textbf{12.38} \\ \hline \hline
\textbf{} & 0.0001/ 1.37 & 1.8e-10 & 16.9 & 29.7 & 1.6e-12 & 10.6 & 2.5 & \textbf{11.87} & 1.0e-15 & 10.6 & 1.2 & \textbf{25.82} \\ \cline{2-13} 
\textbf{5 000 x 200} & 0.0005/ 5.59 & 3.5e-10 & 17.9 & 30.9 & 7.9e-13 & 12.0 & 5.1 & \textbf{6.11} & 2.3e-15 & 12.0 & 2.4 & \textbf{12.86} \\ \cline{2-13} 
\textbf{} & 0.001/ 10.03 & 7.8e-10 & 17.9 & 31.5 & 6.3e-13 & 11.1 & 5.3 & \textbf{5.93} & 4.6e-15 & 11.1 & 2.5 & \textbf{12.64} \\ \hline \hline
\textbf{} & 0.0001/ 1.62 & 5.0e-10 & 11.4 & 64.4 & 9.8e-14 & 13.0 & 11.5 & \textbf{5.58} & 1.9e-15 & 13.0 & 5.0 & \textbf{12.86} \\ \cline{2-13} 
\textbf{10 000 x 300} & 0.0005/ 6.6 & 2.2e-09 & 14.5 & 66.4 & 7.0e-13 & 12.0 & 11.5 & \textbf{5.77} & 9.0e-15 & 12.0 & 5.5 & \textbf{12.07} \\ \cline{2-13} 
\textbf{} & 0.001/ 11.88 & 3.3e-09 & 15.7 & 67.5 & 2.8e-12 & 11.4 & 11.7 & \textbf{5.78} & 1.9e-12 & 11.4 & 5.6 & \textbf{12.02} \\ \hline \hline
\textbf{} & 0.0001/ 4.24 & 6.0e-10 & 20.0 & 267.5 & 3.4e-14 & 13.5 & 42.6 & \textbf{6.29} & 4.3e-15 & 12.99 & 22.5 & \textbf{11.88} \\ \cline{2-13} 
\textbf{10 000 x 3 000} & 0.0005/ 16.58 & 5.7e-10 & 19.0 & 264.2 & 9.0e-14 & 12.3 & 41.6 & \textbf{6.35} & 2.0e-12 & 11.95 & 25.0 & \textbf{10.55} \\ \cline{2-13} 
 & 0.001/ 28.51 & 7.0e-10 & 18.9 & 266.3 & 2.9e-14 & 12.0 & 43.9 & \textbf{6.07} & 5.5e-14 & 11.01 & 26.0 & \textbf{10.23} \\ \hline \hline
 & 0.0001/ 6.18 & 7.0e-09 & 20.1 & 594.2 & 5.6e-14 & 13.4 & 94.2 & \textbf{6.31} & 1.5e-12 & 13.0 & 54.8 & \textbf{10.84} \\ \cline{2-13} 
\textbf{10 000 x 8 000} & 0.0005/ 23.72 & 2.3e-08 & 19.2 & 596.5 & 1.0e-13 & 12.4 & 95.3 & \textbf{6.26} & 2.4e-12 & 12.0 & 61.8 & \textbf{9.66} \\ \cline{2-13} 
 & 0.001/ 39.90 & 2.4e-08 & 18.0 & 584.1 & 4.0e-14 & 12.0 & 98.8 & \textbf{5.91} & 5.1e-14 & 11.0 & 63.1 & \textbf{9.26}
\end{tabular}}}
\end{table}

\begin{table}[]
\centering
\caption{Results for simulations with matrices of different size and the three tested methods in 
C code with mex interface. Error (Err.) and running times are shown for each of them. Speedup with 
respect to LP is shown for SRF and Proposed. Furthermore, we point out that the error of SRF and Proposed 
is three to five orders of magnitude lower than that of LP.}
\label{tab:comps2}
\resizebox{\textwidth}{!}{{
\begin{tabular}{c|r|r|r|r|r|r|r|r|r}
 &  & \multicolumn{2}{c|}{\textbf{LP}~\cite{quattoni2009efficient}} & \multicolumn{3}{c|}{\textbf{SRF}~\cite{chau2018steffensen}} & \multicolumn{3}{c}{\textbf{Proposed}} \\  \hline
\textbf{Matrix Size} & \textbf{$\alpha$ / sp} & \textbf{Err.} & \textbf{Time(s)} & \textbf{Err.} & \textbf{Time(s)} & \textbf{Speedup} & \textbf{Err.} & \textbf{Time(s)} & \textbf{Speedup} \\ \hline  \hline
 & 0.0001/ 1.02 & 4.2e-10 & 0.048 & 2.0e-12 & 0.007 & \textbf{6.67} & 1.9e-16 & 0.005 & \textbf{10.21} \\ \cline{2-10} 
\textbf{2 000 x 100} & 0.0005/ 4.13 & 4.10e-10 & 0.046 & 1.0e-12 & 0.020 & \textbf{2.30} & 7.5e-16 & 0.011 & \textbf{4.17} \\ \cline{2-10} 
\textbf{} & 0.001/ 7.51 & 4.1e-10 & 0.047 & 1.0e-12 & 0.022 & \textbf{2.14} & 1.5e-15 & 0.012 & \textbf{3.94} \\ \hline  \hline
\textbf{} & 0.0001/ 1.37 & 2.6e-08 & 0.25 & 1.7e-12 & 0.055 & \textbf{4.54} & 6.4e-16 & 0.033 & \textbf{7.46} \\ \cline{2-10} 
\textbf{5 000 x 200} & 0.0005/ 5.59 & 2.5e-08 & 0.25 & 5.6e-13 & 0.121 & \textbf{2.04} & 2.3e-15 & 0.066 & \textbf{3.75} \\ \cline{2-10} 
\textbf{} & 0.001/ 10.03 & 2.5e-08 & 0.25 & 6.5e-13 & 0.123 & \textbf{2.02} & 4.5e-15 & 0.067 & \textbf{3.71} \\ \hline  \hline
\textbf{} & 0.0001/ 1.62 & 2.3e-07 & 0.78 & 8.7e-14 & 0.352 & \textbf{2.22} & 1.9e-15 & 0.19 & \textbf{4.04} \\ \cline{2-10} 
\textbf{10 000 x 300} & 0.0005/ 6.6 & 2.3e-07 & 0.78 & 6.1e-14 & 0.386 & \textbf{2.02} & 9.3e-15 & 0.21 & \textbf{3.73} \\ \cline{2-10} 
\textbf{} & 0.001/ 11.88 & 2.3e-07 & 0.78 & 3.0e-12 & 0.393 & \textbf{1.99} & 1.9e-12 & 0.21 & \textbf{3.67} \\ \hline  \hline
\textbf{} & 0.0001/ 4.24 & 2.2e-06 & 8.95 & 3.2e-13 & 4.6356 & \textbf{1.93} & 4.3e-15 & 2.51 & \textbf{3.56} \\ \cline{2-10} 
\textbf{10 000 x 3 000} & 0.0005/ 16.58 & 2.2e-06 & 8.95 & 2.2e-12 & 5.0528 & \textbf{1.77} & 2.1e-12 & 2.72 & \textbf{3.29} \\ \cline{2-10} 
 & 0.001/ 28.51 & 2.2e-06 & 8.93 & 3.9e-13 & 5.1794 & \textbf{1.72} & 5.8e-14 & 2.80 & \textbf{3.19} \\ \hline  \hline
 & 0.0001/ 6.18 & 5.8e-06 & 25.07 & 1.7e-12 & 13.6148 & \textbf{1.84} & 1.5e-12 & 7.35 & \textbf{3.41} \\ \cline{2-10} 
\textbf{10 000 x 8 000} & 0.0005/ 23.72 & 5.8e-06 & 25.08 & 1.9e-12 & 15.0689 & \textbf{1.66} & 2.4e-12 & 8.09 & \textbf{3.10} \\ \cline{2-10} 
 & 0.001/ 39.90 & 5.8e-06 & 25.07 & 3.9e-13 & 15.5143 & \textbf{1.62} & 5.3e-14 & 8.40 & \textbf{2.99}
\end{tabular}}}
\end{table}

For our first set of experiments, we compared pure matlab implementations of our proposed 
method (denoted Proposed) against \cite{sra2011fast},  denoted as general root-finding (GRF), 
and~\cite{chau2018steffensen, chau-l1Infty-codes}, denoted as Steffensen root-finding (SRF). 
Unfortunately, we could not find a public implementation of~\cite{sra2011fast}, and thus, 
we coded our own Matlab version using the \texttt{fzero} function as root search method as suggested in~\cite{sra2011fast}. 

For our second set of comparisons, we tested our proposed method implemented in C 
with a Matlab MEX interface against a similar implementation of SRF~\cite{chau2018steffensen, chau-l1Infty-codes} 
and the linear programming (LP) based method described in Section \ref{sec:quattoni}, which was obtained 
from~\cite{quattoniCode} and also has a MEX interface. None of the C-Mex implementations make use of any 
type of parallelization (SIMD, CUDA, etc.). 

We chose $\text{tol}_c=10^{-12}$ and $\text{tol}_u=10^{-8}$ in~(\ref{eqn:st-tol}) for the {SRF} method. The projections onto the $\ell_1$-ball, needed for the evaluation of the search function in \cite{sra2011fast} and in our algorithm, are implemented using the 
Michelot algorithm~\cite{michelot1986finite}. This algorithm was chosen since it can be implemented efficiently \cite[Section 3.2.2]{rodriguez2017parallel},\cite[Section 3.1]{rodriguez-2018-accelerated} and, for small size projections, we have empirically observed that it has better computational performance than the alternatives mentioned in Section \ref{sec:l1ball}.


We simulated 100 realizations of five different sizes for the matrix $B$, namely\footnotemark[1] 2000 $\times$ 100, 5000 $\times$ 200, 10000 $\times$ 300, 10000 $\times$ 3000 and 10000 $\times$ 8000. For the constraint parameter, we took $\tau=\alpha \| B \|_{\infty,1}$. We considered\footnotemark[1] $\alpha \in  \{10^{-4},\, 5 \times 10^{-4},\, 10^{-3}\}$,  to experimentally obtain approximate sparsity percentages (percentage of non-zero rows) of 1, 5 and 10\%, respectively.

\begin{figure*}
  \centering
      \includegraphics[trim={2cm 14cm 5cm 0},clip, width=1\textwidth]{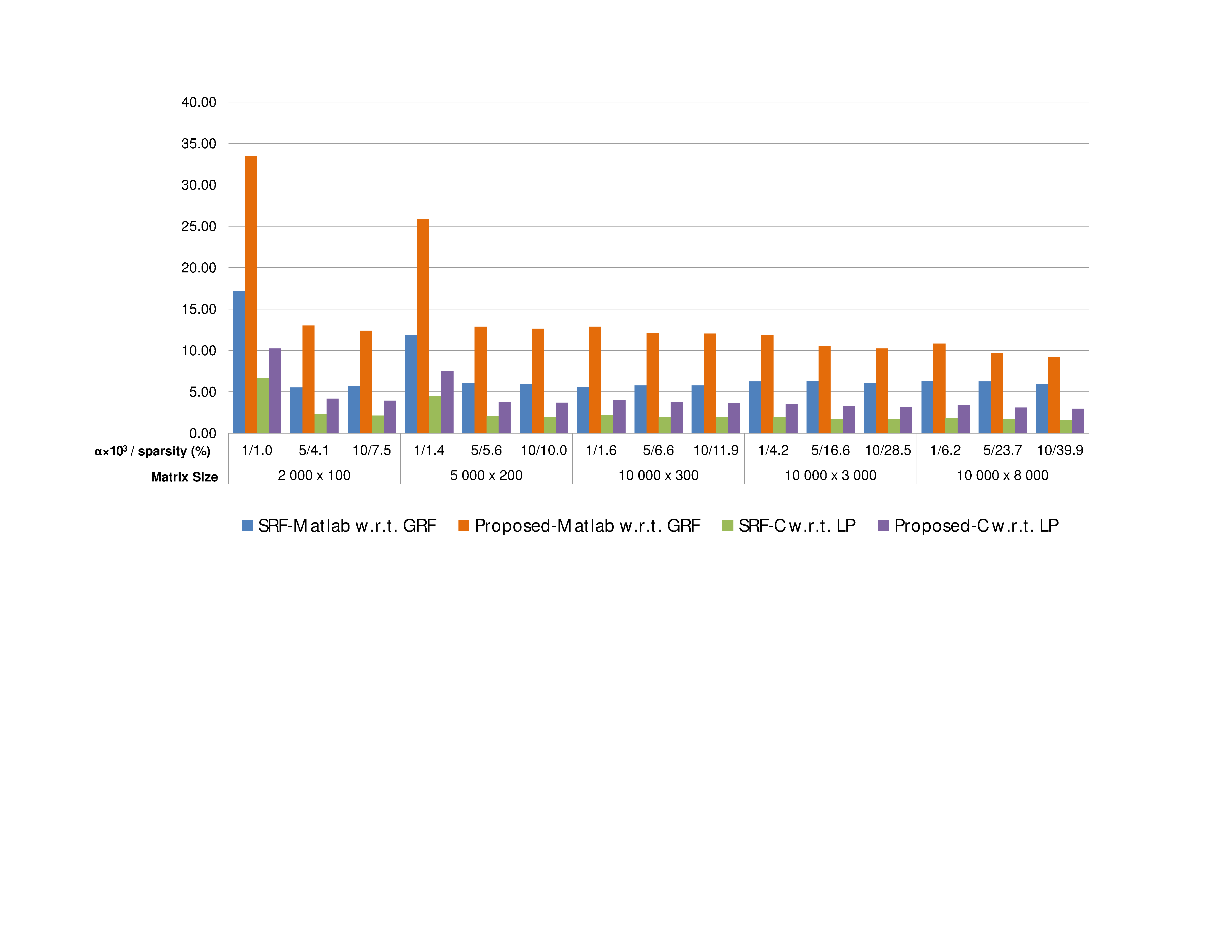}
      \caption{{Bar plot showing the average speedup with respect to GRF for SRF-Matlab (blue), proposed-Matlab (orange), and with respect to LP for SRF-C (green) and Proposed-C (purple) for the different matrix size and $\alpha$/sparsity values. See also Tables~\ref{tab:comps} and \ref{tab:comps2}.}}
      \label{fig:speedup}
\end{figure*}

As discussed in~\cite{sra-2012-fast}, if $\| B\|_{\infty,1} > \tau$ then, at the optimum, the inequality 
constraint of~\eqref{eq:proj_1inf_ball1} is active. For all our test, we made sure that 
$\| B \|_{\infty,1} > \tau$, and thus we measured the error of the solution as $| \, \|X\|_{\infty,1} -\tau |$. 
We additionally computed the number of iterations, execution time, and a sparsity value as the percentage of non-zero rows. 
As all methods arrive to the same value of sparsity, only a single value is shown for each $\alpha$. 
The results averaged over the 100 realizations for the different matrix dimensions are shown in 
Table~\ref{tab:comps} and the average speedups are shown in Figure~\ref{fig:speedup}.

\subsection{Comparisons on fMRI LASSO application}
\label{sec:fmri}

\begin{figure*}
  \centering
      \includegraphics[width=\textwidth]{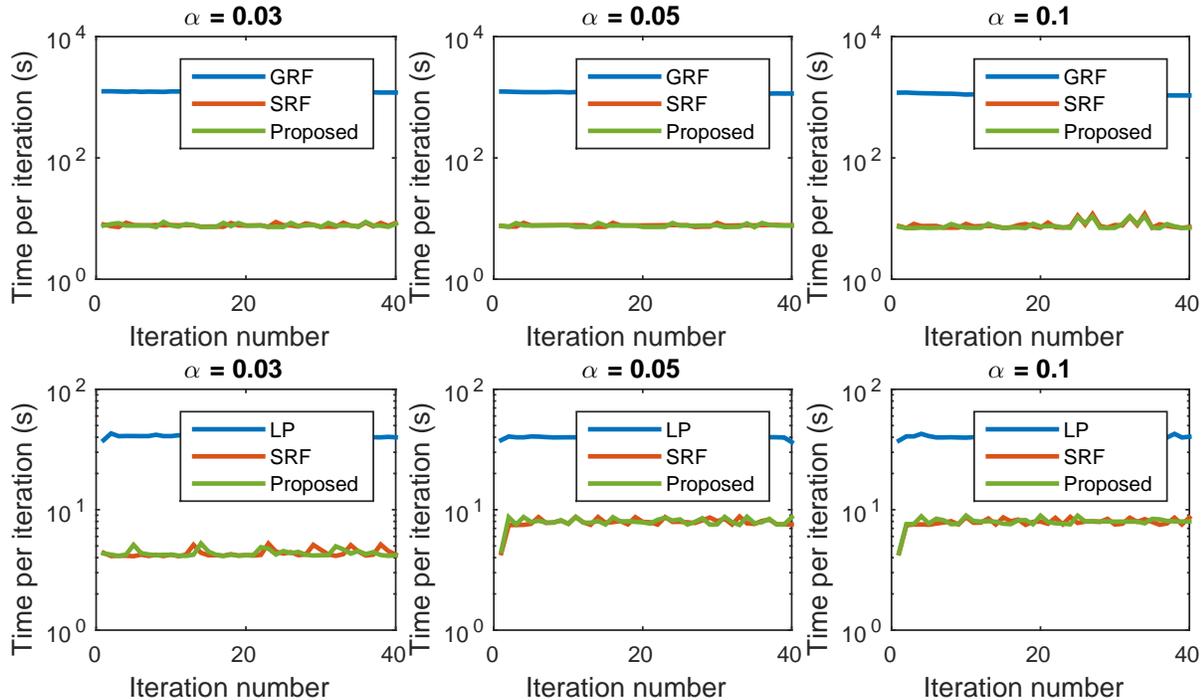}
      \caption{{Time (in seconds) per PGD iteration for each of the methods for solving 
      the $\ell_{\infty,1}$ projection problem associated with the fMRI experiment for Matlab coded methods (Top row) and C coded methods (Bottom row).}}
      \label{fig:comp_fmri}
\end{figure*}

\begin{table}[]
\centering
{
\caption{Mean time in seconds and speedup with respect to GRF for the fMRI experiment (Matlab coded methods).}
\label{tab:comp_fmri-1}
\begin{tabular}{c|c||c|c||c|c}
 & \textbf{GRF}~\cite{sra2011fast} & \multicolumn{2}{c||}{\textbf{SRF}~\cite{chau2018steffensen}} & \multicolumn{2}{c}{\textbf{Proposed}} \\ \hline
\textbf{$\alpha$ / sp} & \textbf{Time(s)} & \textbf{Time(s)} & \textbf{speedup} & \textbf{Time(s)} & \textbf{speedup} \\ \hline
\textbf{0.03/12.5\%} & 1220.7 & 7.78 & \textbf{156.9} & 7.75 & \textbf{157.51} \\ \hline
\textbf{0.05/31.4\%} & 1189.3 & 7.71 & \textbf{154.25} & 7.66 & \textbf{155.26} \\ \hline
\textbf{0.1/ 58.1\%} & 1109.5 & 7.89 & \textbf{140.62} & 7.68 & \textbf{144.47}
\end{tabular}}
\end{table}

\begin{table}[]
\centering
{
\caption{Mean time in seconds and speedup with respect to LP for the fMRI experiment (C coded methods).}
\label{tab:comp_fmri-2}
\begin{tabular}{c|c||c|c||c|c}
 & \textbf{LP}~\cite{quattoni2009efficient} & \multicolumn{2}{c||}{\textbf{SRF}~\cite{chau2018steffensen}} & \multicolumn{2}{c}{\textbf{Proposed}} \\ \hline
\textbf{$\alpha$ / sp} & \textbf{Time(s)} & \textbf{Time(s)} & \textbf{speedup} & \textbf{Time(s)} & \textbf{speedup} \\ \hline
\textbf{0.03/12.5\%} & 40.3 & 4.36 & \textbf{9.24} & 4.36 & \textbf{9.24} \\ \hline
\textbf{0.05/31.4\%} & 39.8 & 7.86 & \textbf{5.06} & 7.94 & \textbf{5.01} \\ \hline
\textbf{0.1/ 58.1\%} & 39.9 & 7.86 & \textbf{5.08} & 7.92 & \textbf{5.04} 
\end{tabular}}
\end{table}

We tested the computational improvements of our algorithm in the cognitive task described in~\cite{mitchell2008predicting}. The applications consists of predicting the neural functional magnetic resonance image (fMRI) response associated to a particular word based on co-occurrence features of this word with a dictionary of words whose response is already know. 
In ~\cite{mitchell2008predicting}, the co-occurrence with a hand-crafted set of 25 verbs was used as features in the prediction problem, whereas ~\cite{liu2009blockwise} showed improvements by using a larger dictionary and MTL to select the best features.

For our tests, we selected the 18 noun words with their corresponding fMRI images and used co-occurrence values of these with a dictionary of 10000 words gathered from Wikipedia and BBC, which was obtained from~\cite{dissect}. These set of 18 words was selected because the other words of the dataset were not present in the corpus used for the co-occurrence matrix calculation.

We subsampled the fMRI images by half resulting in approximate problem dimension of $K=10000$, $N=18$ and $M=10000$. $\alpha$ values of 0.03, 0.05 and 0.1 were considered in order to obtain sparsity values of 12.5\%, 31.4\% and 58.1\%, respectively.

The improvements that can be obtained with the use of mixed norms in MTL has already been demonstrated~\cite{liu2009blockwise}, 
so we focus only on computational metrics. We compare the time needed to solve~\eqref{eq:mtl} by using GRF, SRF and the 
proposed method as the projection operators in algorithm~\ref{al:pgd}. For PGD, we use the minConf Matlab 
library~\cite{schmidt2009optimizing,minConf} which uses an Armijo inexact line search~\cite{NoceWrig06} 
for choosing the step size. We consider a fixed maximum of 40 iterations (due to the long computational times). 
The results for each case ($\alpha = \{ 0.03, 0.05, 0.1\}$) are shown in Figure~\ref{fig:comp_fmri} for our proposed method
as well as for all other (GRF, SRF and LP) methods; moreover, in Tables~\ref{tab:comp_fmri-1} and~\ref{tab:comp_fmri-2} 
we list the corresponding mean time and speedup for the MATLAB and C coded implementations respectively.
%

\subsection{Discussion}
\label{sec:discussion}

As can be observed from the results of Section~\ref{sec:initial}, the initial point $\gamma_0$ has impact only at low $\alpha$ (and therefore $\tau$) values. This is easily explainable, as at high $\tau$ values the solution of~\eqref{eq:shrink} is zero. Accordingly, as suggested in Section~\ref{sec:algoritmo}, it is better to first evaluate the conditions on $\max_{n,k} \left\lbrace b_{n,k} \right\rbrace$ and $\|\mathbf{b}_i\|_{\infty}$ to avoid unnecessary shrinkage operations, as these comparisons do not incur a great computational cost. When the initial point is different from zero, we see that the number of iterations and the computational time are slightly reduced. For our tests, when $\gamma_0$ is not zero, an average reductions of 2.2 iterations ($-15\%$) and of $65\%$ in time are achieved.

As indicated by the results of Section~\ref{sec:simulations}, the proposed Newton method tends to require fewer iterations than GRF~\cite{sra2011fast}, {although
this is not directly comparable as the stopping criteria of the
fzero function used in GRF (only implemented in Matlab) is different
from the one we are using in ~\ref{al:proposed} (implemented in Matlab and C-MEX)} and because our proposed algorithms are obtaining smaller errors. Additionally, the Newton-based method performs the same number of iterations as the SRF method \cite{chau2018steffensen}, but the computational time is substantially lower. This is explained by the fact that in the Newton-based method, only one function evaluation is performed compared to the two evaluations performed for the SRF method. Overall, the proposed Newton-based algorithm obtains speedups of 8 or more with respect to GRF (compared to the average of 5 obtained with SRF). Higher speedups are obtained at low $\alpha$ (sparser) values, as this is where the initial point guess is more effective~\cite{chau2018steffensen}. For those cases, the speedup of our proposed method goes up to 14 $\sim$ 25. {Likewise, in the C code comparisons, our proposed method obtained speedups of 3 $\sim$ 10 with respect to that of LP. The same effect of the initial point guess providing higher speedups at lower $\alpha$ values is also observed.} 

\begin{figure}
  \centering
      \includegraphics[width=0.6\textwidth]{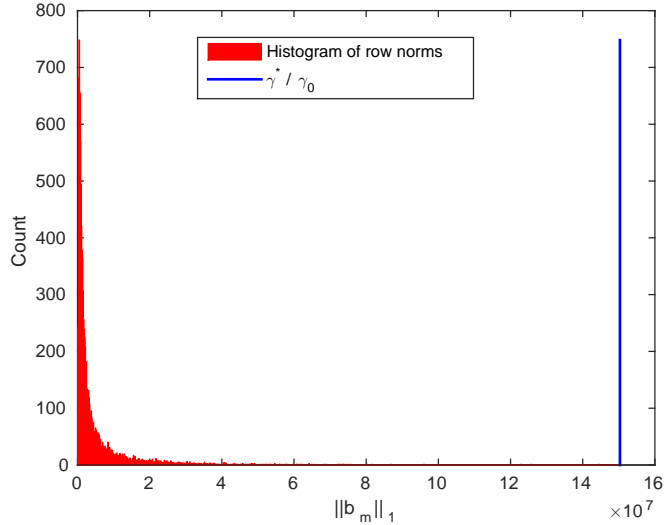}
      \caption{Distribution of $||\mathbf{b_m}||_1$ values (red) and optimal $\gamma$ value (blue).}
      \label{fig:explain_fmri}
\end{figure}

Regarding the fMRI experiments, we observe that both our previous SRF method and our novel Newton-based methods obtain considerable speedups of 120$\sim$130. This speedup reduced the total computational time for the whole PGD method from around 10 hours (GRF) to approximately 3 minutes (SRF and Proposed). {In the C-code comparisons, we obtained speedups of 5$\sim$9 with respect to the LP method.} The difference in performance with respect to the simulations is explained by the data distribution. As observed in Figure~\ref{fig:explain_fmri}, the $\ell_1$-norms of the rows of the data follow a Laplacian-like distribution, while the simulations considered uniformly distributed data.  Accordingly, in the fMRI data, the pruning strategy mentioned in Section~\ref{sec:search_fun} is able to greatly reduce the number of rows to be processed. Furthermore, the initial guess $\gamma_0$ obtained with the SRF and proposed methods is very close to the optimal $\gamma^*$, and so we noticed that only around two root-search iterations were needed for each projection. Although for this particular dataset the proposed method performs comparably to SRF, the simulations showed that for setting where the data is more uniformly distributed the difference in computation can be considerable.

\section{Conclusion}
\label{sec:conclusions}

We have presented a new algorithm for efficient projection onto the $\ell_{\infty,1}$-norm ball which exploits the particular structure of this problem to improve over previous methods. The algorithm, based on a Newton root-finding approach, capitalizes on an approximation of the derivative of the  search function obtained from Moreau's decomposition and duality theory. Our proposed method obtains speedups of eight or more with respect to previous state-of-the-art methods, while achieving smaller errors in our simulations.
When we applied our proposed algorithm to a multi-task Lasso (MTL) using real
fMRI data, we obtained considerable speedups. Furthermore, the MTL test also highlights the impact 
of two key aspects of our proposed algorithm, namely our initial guess and the pruning steps. 
When the distribution of the data is favorable, these aspects have a very significant positive impact on the 
overall computational performance of our proposed algorithm.
{Implementations of the proposed algorithms in both Matlab and C are provided in~\cite{chau-l1Infty-codes}, and 
a Python implementation will be included in a future release of the SPORCO library~\cite{sporco}. The data 
and scripts necessary for reproducing the experiments reported here are also available in~\cite{chau-l1Infty-codes}.}

\section*{Acknowledgment}

This research was supported by the ``Programa Nacional de Innovaci\'on para la Competitividad y Productividad'' (Inn\'ovate Per\'u) Program, 169-Fondecyt-2015, and by the U.S. Department of Energy through the LANL/LDRD Program.





%

\bibliographystyle{IEEEtran}
\bibliography{mixed_norm}

\end{document}